\documentclass[12pt,reqno]{amsproc}
\usepackage[utf8]{inputenc}
\usepackage[margin=.9in]{geometry}
\usepackage{amsmath,mathrsfs}
\usepackage{graphicx}
\usepackage[margin=10pt,font=small,labelfont=bf]{caption}
\usepackage[all,2cell,cmtip]{xy}
\UseAllTwocells

\usepackage[sc]{mathpazo}\renewcommand{\mathbf}{\mathbold}

\newcommand{\set}[1]{\left\{#1\right\}}

\newcommand{\sepline}{}

\usepackage{hyperref}
\hypersetup
{
	colorlinks,	%
	citecolor=red,%
	filecolor=blue,%
	linkcolor=blue,%
	urlcolor=blue
}

\usepackage{amssymb}
\usepackage{xcolor}

\newtheorem{theorem}{Theorem}[section]

\theoremstyle{definition}
\newtheorem{definition}[theorem]{Definition}

\newcommand{\Barc}{\text{\bf Bar}}
\newcommand{\Z}{\mathbb{Z}}
\newcommand{\R}{\mathbb{R}}

\title[Vectorization Methods in Topological Data Analysis]{A Survey of Vectorization Methods\\ in Topological Data Analysis}

\author[Ali]{Dashti Ali}
\author[Asaad]{Aras Asaad}
\author[Jimenez]{Maria Jose Jimenez}
\author[Nanda]{Vidit Nanda}\thanks{Corresponding Author: nanda@maths.ox.ac.uk}
\author[Paluzo-Hidalgo]{Eduardo Paluzo-Hidalgo}
\author[Soriano-Trigueros]{Manuel Soriano-Trigueros}\thanks{Authors listed in alphabetical order}

\begin{document}

\maketitle

\begin{abstract}
    Attempts to incorporate topological information in supervised learning tasks have resulted in the creation of several techniques for vectorizing persistent homology barcodes. In this paper, we study thirteen such methods. Besides describing an organizational framework for these methods, we comprehensively benchmark them against three well-known classification tasks. Surprisingly, we discover that the best-performing method is a simple vectorization, which consists only of a few elementary summary statistics. Finally, we provide a convenient web application which has been designed to facilitate exploration and experimentation with various vectorization methods.
\end{abstract}

\section*{Introduction}

\begin{figure}
    \includegraphics[width=.9\linewidth]{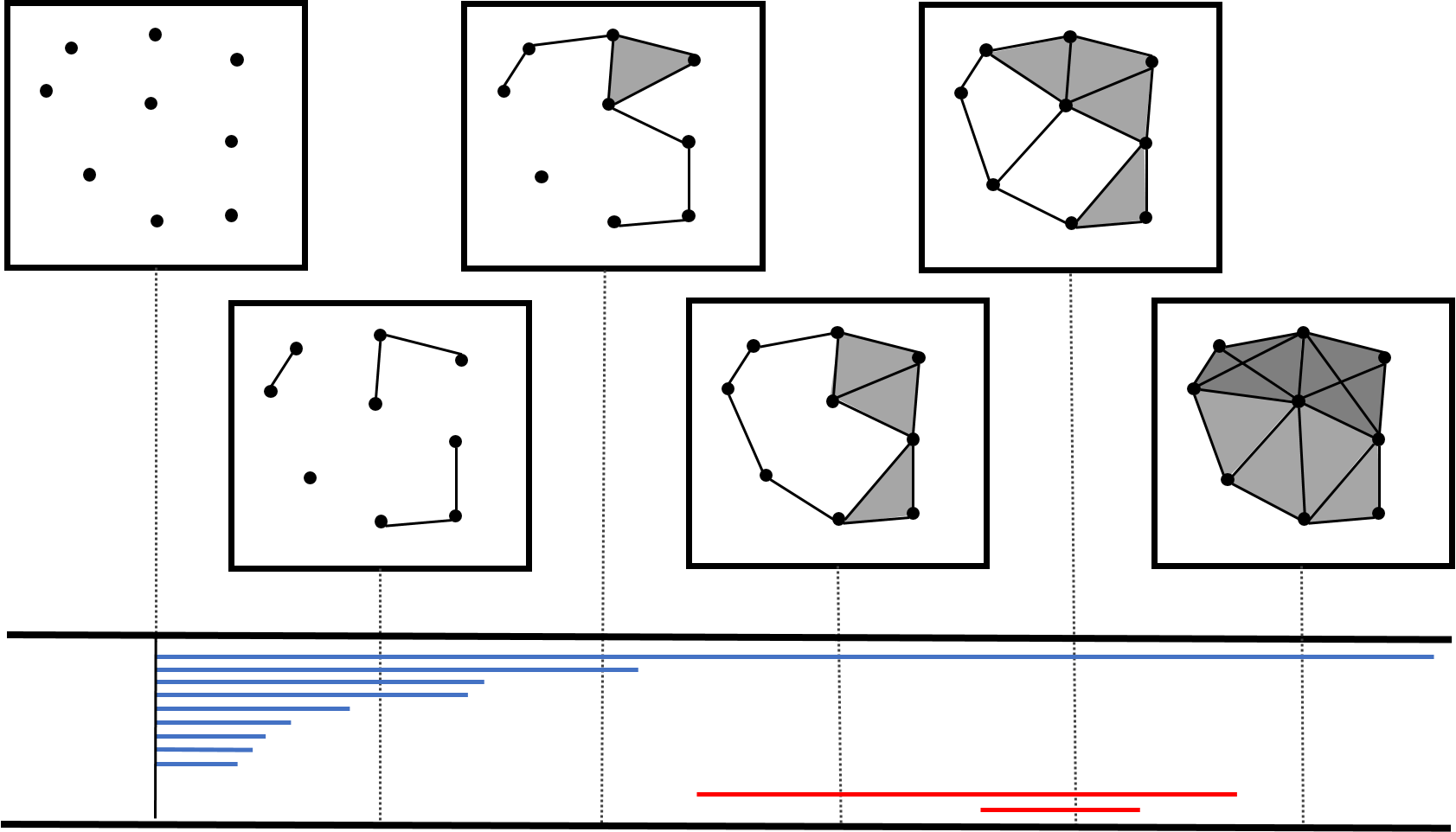}
    \caption{An increasing family of cell complexes built around a point cloud dataset; the associated barcode in dimensions 0 (blue) and 1 (red) catalogues the connected components and cycles respectively.}\label{fig:barcode}
\end{figure}

Propelled by deep theoretical foundations and a host of computational breakthroughs, topological data analysis emerged roughly three decades ago as a promising method for extracting insights from unstructured data \cite{ghrist, carlsson:09, nansaz, oudot}. The principal instrument of the enterprise is {persistent homology}; this consists of three basic steps, each relying on a different branch of mathematics. 
\begin{enumerate}
    \item {\em Metric geometry}: construct an increasing family $\set{X_t}$ of cell complexes around the input dataset $X$, where the indexing $t$ is a scale parameter in $\R_{\geq 0}$.
    \item {\em Algebraic topology}: compute the $d$-th homology vector spaces $H_d(X_t)$ for scales $t$ in $\R_{\geq 0}$ and dimensions $d$ in $\Z_{\geq 0}$.
    \item {\em Representation theory}: decompose each family of vector spaces $\set{H_d(X_t) \mid t \geq 0}$ into irreducible summands, thus producing a {\bf barcode}. 
\end{enumerate}
The resulting barcodes are finite multisets of real intervals $[p,q] \subset \R$, which admit concrete geometric interpretations in low dimensions --- see Figure \ref{fig:barcode}. The ultimate goal is to infer the coarse geometry of $X$ across various scales by examining the longer intervals in its barcodes. Crucially, once the method for constructing $\set{X_\epsilon}$ from $X$ has been fixed, the entire persistent homology pipeline is unsupervised: one requires neither labelled data nor hyperparameter tuning to produce barcodes from $X$. 

At the other end of the data analysis spectrum lies supervised machine learning using contemporary neural networks, which are replete with billions of tunable parameters and gargantuan training datasets \cite{neuralbook}. The practical aspects of deep neural networks appear to be light years ahead of the underlying theory. It nevertheless remains the case that machine learning has driven astonishing progress in the systematic automation of several important classification tasks. One direct consequence of these success stories is the irresistible urge to combine topological methods with machine learning. The most common avenue for doing so is to turn barcodes into vectors (lying in a convenient Euclidean space) which then become input for suitably-trained neural networks. 

The good news, at least from an engineering perspective, is that barcodes are inherently combinatorial objects, and as such, they are remarkably easy to vectorize. Several dozen vectorization methods have been proposed across the last decade, and new ones continue to appear with alarming frequency and increasing complexity --- the reader will encounter thirteen of them here. The bad news, on the other hand, comes in the form of three serious challenges which must be confronted by those who build or use such vectorizations:
\begin{enumerate}
    \item Given the large number of options, even established practitioners are not aware of all the vectorization techniques; similarly, knowledge of which vectorizations are suitable for which types of data is difficult -- if not impossible -- to glean from the published literature.
    \item There is a natural metric between barcodes called the {\em bottleneck distance}; when it is endowed with this metric, the space of barcodes becomes infinite-dimensional and highly nonlinear. As such, it does not admit any faithful embeddings into finite-dimensional vector spaces.
    \item Even the {\em stable} vectorizations, which preserve distances by mapping barcodes into infinite-dimensional vector spaces, may suffer from a lack of discriminative power in practice: by design, they are poor at distinguishing between datasets whose coarse structures are similar and whose differences reside in finer scales.  
\end{enumerate}

\subsection*{In This Paper} Here we seek to comprehensively describe, catalogue and benchmark vectorization methods for persistent homology barcodes. The first contribution of this paper is the following taxonomy of the known methods, which we hope will serve as a convenient organizational framework for beginners and experts alike ---
\begin{enumerate}
    \item {\em Statistical} vectorizations: these summaries consist of basic statistical quantities;
    \item {\em Algebraic} vectorizations: these are generated from polynomials;
    \item {\em Curve} vectorizations: these come from maps $\R \to H$, where $H$ is a vector space;
    \item {\em Functional} vectorizations: these are maps of the form $X \to H$ for $X \neq \R$;
    \item {\em Ensemble} vectorizations: these are generated from collections of training barcodes.
\end{enumerate}
There are unavoidable overlaps between these five categories. When such an overlap occurs, we have placed the given vectorization technique in the earliest relevant category among those in the list above; thus, an algebraic vectorization given by polynomial functions of basic statistical quantities will be placed in category (1) rather than category (2). The reader might claim, quite reasonably, that category (3) should be subsumed into category (4). However, the sheer number of curve-based vectorizations compelled us to set them apart. 

The second contribution of this paper is a comprehensive benchmarking of thirteen vectorization techniques across these five categories on three well-known image classification datasets. These datasets were selected to simultaneously ({\bf a}) provide an increasing level of difficulty for topological methods, and ({\bf b}) to be instantly recognizable for the broader machine learning community. These are: the {\bf Outex} texture database \cite{ojala2002outex}, the {\bf SHREC14} shape retrieval dataset \cite{Pickup2014}, and the {\bf Fashion-MNIST} database \cite{xiao2017fashion}. Surprisingly, the best-performing vectorization {\em in all three cases} is a rather na{\"i}ve one obtained by collecting basic statistical quantities associated to (the multiset of) intervals in a given barcode.

Our third contribution is a companion web application which computes and visualizes all thirteen vectorization techniques which have been investigated in this paper. In addition to running online\footnote{\url{https://persistent-homology.streamlit.app}}, this web app can also be downloaded\footnote{\url{https://github.com/dashtiali/vectorisation-app}} and run locally on more challenging datasets. 

\subsection*{Not In This Paper} Vectorization methods form but a small part of the ever expanding interface between topological data analysis and machine learning. As such, there are several related techniques which are not benchmarked here. The precise inclusion criteria for our study in this paper are as follows.
\begin{enumerate}
    \item We restrict our attention to those methods which produce genuine vectors from barcodes. In particular, kernel methods \cite{reininghaus2015, carriere2017sliced} are beyond the scope of this paper.
    \item We only consider those vectorizations that are either straightforward for us to implement, or have an easily accessible and trusted implementation. For instance, path signature based vectorizations \cite{Nanda2020, giusti2021} are excluded.
    \item We do not compare machine learning architectures designed for the explicit purpose of inferring (persistent) homology \cite{carriere2020perslay, hofer2019, keros22}. 
    \item We do not touch upon various attempts to design or study neural networks using tools from topological data analysis \cite{moor2020,carlsson2020topological}.
    \item Finally, even among methods which satisfy the first four criteria, we have discarded techniques which regularly obtained a classification accuracy below fifty percent.
\end{enumerate}

\subsection*{Similar Efforts}  The authors of \cite{pun2018} have summarised -- but not compared -- several vectorization and kernel methods for barcodes. Another summary (sans comparison) may be found in \cite{som2020}, with emphasis on metric aspects of the chosen vectorizations. The work of \cite{Lawson2022} describes a common overarching framework for what we have called curve vectorizations here. More recently, \cite{Perea2021} and \cite{conti2022} have described and compared five and four vectorization methods respectively. 

\subsection*{Outline} Notation and preliminaries involving barcodes are established in Section \ref{sec:barcodes}. In Sections \ref{sec:vecs} and \ref{sec:datasets} we introduce the thirteen vectorizations (suitably organised into our taxonomy) and the three datasets. Section \ref{sec:results} contains the results of our experiments whose finer details have been relegated to Appendices \ref{sec:impl} and \ref{app:shrec}. We provide a description of the web app in Section \ref{sec:webapp} and some brief concluding remarks in Section \ref{sec:conclusion}.

\section{Persistence Barcodes from Data}\label{sec:barcodes}

At its core, persistent homology studies sequences of finite-dimensional vector spaces $V=\set{V_i \mid 0 \leq i \leq n}$ and linear maps $a=\set{a_i:V_{i-1} \to V_{i} \mid 1 \leq i \leq n}$:
\[
\xymatrixcolsep{.7in}
\xymatrix{
V_0 \ar@{->}[r]^{a_1} & V_1 \ar@{->}[r]^{a_2} & \cdots  \ar@{->}[r]^{a_n} & V_n.
}
\]
Such sequences $(V,a)$ are called {\em persistence modules}. Among the simplest examples are {\em interval} modules --- for each pair of integers $p \leq q$ with $[p,q] \subset [0,n]$, the corresponding interval module $(I^{[p,q]},c^{[p,q]})$ has 
\[
\dim I^{[p,q]}_i=\begin{cases} 1 & \text{if } p \leq i \leq q \\
0 & \text{otherwise;} \end{cases}
\]
similarly, the map $c^{[p,q]}_i$ is the identity whenever $p+1 \leq i \leq q$ and zero otherwise.

\subsection{Structure and Stability}
Every persistence module decomposes into a direct sum of interval modules. In particular, we have the following {structure theorem} \cite{zomcar, cdsgo}.

\begin{theorem}
For every persistence module $(V,a)$, there exists a unique set $\Barc(V,a)$ of subintervals of $[0,n]$ along with a unique function $\Barc(V,a) \to \Z_{>0}$ denoted $[p,q] \mapsto \mu_{p,q}$ for which we have an isomorphism
\[
(V,a) \simeq \bigoplus_{[p,q] \in \Barc(V,a)} \left(I^{[p,q]},c^{[p,q]}\right)^{\mu_{p,q}}.
\]
\end{theorem}
\noindent Thus, the algebraic object $(V,a)$ may be fully recovered (up to isomorphism) from purely combinatorial data consisting of the set of intervals $\Barc(V,a)$ and the multiplicity function $\mu$. Alternately, one may view $\Barc(V,a)$ itself as a multiset with $\mu_{p,q}$ copies of each interval $[p,q]$. This multiset is called the {\bf barcode} of $(V,a)$. It is often useful in applications to let the vector spaces $V_i$ be indexed by real numbers rather than integers. With this modification in place, $\Barc(V)$ becomes a collection of real intervals $[p,q] \subset \R$. 

The most important property of persistence modules, beyond the structure theorem, is their {\bf stability} \cite{cdsgo}. There is a natural metric on the set of persistence modules called the {\em interleaving distance} and a metric on the set of barcodes called the {\em bottleneck distance}

\begin{theorem}\label{thm:stability}
 The assignment $(V,a) \mapsto \Barc(V,a)$ is an isometry from the space of persistence modules (with interleaving distance) to the space of barcodes (with bottleneck distance).
\end{theorem}

The advantage of this theorem is that barcodes remain robust to (certain types of) perturbations of the original dataset, thus conferring upon the topological data analysis pipeline a degree of noise-tolerance. The significant difficulty from a statistical perspective, however, is that the metric space of persistence barcodes with bottleneck distance is nonlinear --- even averages can not be defined for arbitrary collections of barcodes \cite{frechetmeans, persgeom, higherint}. 

\subsection{Barcodes from Data} Persistence modules arise naturally from a wide class of datasets. The first step in topological data analysis involves imposing the structure of a filtered cell complex -- either simplicial \cite[Chapter 8]{simphom} or cubical \cite{Kaczynski2004} -- from the data \cite{ghrist, carlsson:09, nansaz}. The two most prominent examples of filtered cell complex structures arising from data are as follows.
\begin{enumerate}
    \item Given a finite point cloud $X \subset \mathbb{R}^n$, one constructs a family of increasing simplicial complexes $\set{S_\epsilon \mid \epsilon \geq 0}$ defined as follows. A collection $\set{x_0,\ldots,x_k}$ forms a $k$-simplex in $S_\epsilon$ if and only if the (Euclidean) distance between $x_i$ and $x_j$ is no larger than $\epsilon$ for all $i,j$ in $\set{0,\ldots,k}$. Since there are only finitely many $\epsilon$ values at which new simplices are introduced, the filtration is indexed by a subset of the natural numbers. The collection $S_\epsilon$ is called the {\bf Vietoris-Rips} filtration of $X$. These filtrations can be defined for any metric space in a similar fashion.   
    \item Consider a grayscale image $I$, given in terms of $m \times n$ pixels with intensity values in the set $\set{0,1,\ldots,255}$. This naturally forms a two-dimensional {cubical complex}, which can be endowed with the {\bf upper-star} filtration by intensity values. In particular, each elementary cube of dimension $< 2$ appears at the smallest intensity encountered among the 2-dimensional cubes in its immediate neighbourhood. Higher-dimensional cubical filtrations may be similarly generated from higher-dimensional pixel grids.
\end{enumerate}

Once the given dataset has been suitably modeled by a filtered cell complex, persistence modules are obtained by computing {\bf homology} groups with coeffiecients in a field. The reader who is interested in the definition and computation of homology is urged to either consult standard algebraic topology references such as \cite[Ch 2]{hatcher} or see the more recent \cite{oudot, Edelsbrunner2010, nansaz}. 

A substantial difficulty in topological data analysis is that although persistent homology barcodes can be readily associated with a large class of datasets, the space of all such barcodes is notoriously unpleasant to encounter from a statistical perspective. Fortunately, barcodes are combinatorial objects which can be mapped to Hilbert spaces in a plethora of reasonable ways. Indeed, across the last decade, such {\bf vectorization methods} have been proposed by various authors, and our main purpose in this work is to benchmark many of these methods against standard classification tasks. 

\section{Vectorization Methods for Barcodes}\label{sec:vecs}

Throughout this section, we assume knowledge of the barcode $B:= \Barc(V,a)$ of an $\R$-indexed persistence module along with its multiplicity function $\mu:B \to \Z_{>0}$. We note that for each interval $[p,q]$ in $B$ the numbers $p$ and $q$ are called its {\em birth} and {\em death} respectively, and the length $q-p$ is called its {\em lifespan}.

\subsection{Statistical Vectorizations}

The first and simplest category of vectorizations considered in this paper are generated from basic statistical quantities associated to the given barcode. Variants of the following vectorization have been defined and used on several occasions --- see for instance \cite[sec 2.3]{AADA} , \cite[Sec 6.2.1]{Lawson2022}  and \cite[Sec 4.1.1]{pun2018}.

\begin{definition}\label{def:pers_stats}
The {\bf persistence statistics} vector of $\mu:B \to \Z_{>0}$ consists of: 
\begin{enumerate} 
\item the mean, the standard deviation, the median, the interquartile range, the full range, the $10^\text{th}$, $25^\text{th}$, $75^\text{th}$ and $90^\text{th}$ percentiles of the births $p$, the deaths $q$, the midpoints $\frac{p+q}{2}$ and the lifespans $q-p$ for all intervals $[p,q]$ in $B$ counted with multiplicity;
\item the total number of bars (again counted with multiplicity), and
\item the {\em entropy} of $\mu$, defined as the real number
\[
E_\mu := -\sum_{[p,q] \in B} \mu_{p,q} \cdot \left(\frac{q-p}{L_\mu}\right) \cdot \log \left(\frac{q-p}{L_\mu}\right),
\]
where $L_\mu$ is the weighted sum 
\begin{align}\label{eq:L}
L_\mu := \sum_{[p,q] \in B} \mu_{p,q} \cdot (q-p).
\end{align}
\end{enumerate}
\end{definition}

\noindent\sepline

The entropy from Definition \ref{def:pers_stats}(3) was introduced in \cite{Chintakunta2015, Rucco2016}. Our second statistical vectorization is from \cite{atienza2020}, where entropy has been upgraded to a real-valued piecewise constant function rather than a single number. 

\begin{definition}\label{def:entropy_summary}
The {\bf entropy summary} function of $\mu:B \to \Z_{>0}$ is the map $S_\mu:\R \to \R$ given by
\[
S_\mu(t)=-\sum_{[p,q] \in B} \mathbb{1}_{p \leq t < q} \cdot \mu_{p,q} \cdot \left(\frac{q-p}{L_\mu}\right) \cdot \log \left(\frac{q-p}{L_\mu}\right).
\]
Here $\mathbb{1}_\bullet$ is the indicator function --- it equals $1$ when the conditional $\bullet$ is true and it equals $0$ otherwise. The number $L_\mu$ appearing in the expression above is defined in \eqref{eq:L}.
\end{definition}
\noindent The entropy summary function has also been called the {\em life entropy curve}, e.g., in \cite{Lawson2022}.

\subsection{Algebraic Vectorizations}

The vectorizations in this category are generated using polynomial maps constructed from the barcode $\mu:B \to \Z_{>0}$. 

\noindent\sepline

The first example considered here is from \cite{Adcock2016}. It becomes convenient, for the purpose of defining it, to arbitrarily order the intervals in $B$ as $\set{[p_i,q_i] \mid 1 \leq i \leq n}$ with the understanding that each $[p,q]$ occurs $\mu_{p,q}$ times in this ordered list.

\begin{definition}\label{def:algfun} The ring of {\bf algebraic functions} on $\mu:B \to \Z_{>0}$ consists of all those $\R$-polynomials $f$ in variables $\set{x_1,y_1,\ldots,x_n,y_n}$ for which the following property holds: there exist  polynomials $\set{g_i \mid 1 \leq i \leq n}$ satisfying 
\[
\frac{\partial f}{\partial x_i} + \frac{\partial f}{\partial y_i}=(x_i - y_i) \cdot g_i.
\]
(Here $\partial f/\partial x_i$ indicates the partial derivative of $f$ with respect to $x_i$, and so forth).
\end{definition}

The desired vectorization is obtained by selecting finitely many algebraic functions from this ring and evaluating them at $x_i=p_i$ and $y_i=q_i$ for all $i$. The feature maps generated by making such choices are sometimes called {\em Adcock-Carlsson coordinates} --- see for instance \cite{Perea2022}. Letting $q_\text{max}$ be the maximum death-value encountered among the intervals in $B$, four of the most widely-used algebraic functions are:
\begin{alignat*}{5}
    &f_1=\sum_i p_i(q_i - p_i) \qquad && 
    f_2 =\sum_i(q_\text{max}-q_i)\,(q_i -p_i) \\ 
    &f_3=\sum_i p_i^2(q_i - p_i)^4 \qquad &&
    f_4=\sum_i(q_\text{max}-q_i)^2\,(q_i - p_i)^4
\end{alignat*}

Small changes in the barcode (in terms of bottleneck distance) are liable to create large fluctuations in the associated algebraic functions. The methods of tropical geometry were used in \cite{Kalisnik2019} to address the bottleneck instability of algebraic functions. In this setting, the standard polynomial operations $(+,\times)$ are systematically replaced by $(\max,+)$. To define the resulting vectorization, we once again use an ordering $\set{[p_i,q_i] \mid 1 \leq i \leq n}$ of the intervals in $B$.

\begin{definition}\label{def:tropical} A {\bf tropical coordinate function} for $\mu:B \to \Z_{>0}$ is a function $F$ of variables $\set{x_1,y_1,\ldots,x_n,y_n}$ which is both tropical and symmetric as described below.
\begin{enumerate}
    \item {\em Tropical}: there is an expression for $F$ which uses only the operations max, min, $+$ and $-$ on the variables $\set{x_i}$ and $\set{y_i}$. 
    \item {\em Symmetric}: any permutation of $\set{1,\ldots,n}$, when applied to both $\set{x_i}$ and $\set{y_i}$, leaves $F$ unchanged.
\end{enumerate}
\end{definition}
Let $\lambda_i$ be the lifespan $q_i-p_i$ of the $i$-th interval in $B$. To generate feature maps from the tropical coordinate functions described above, one simply evaluates them at $x_i=\lambda_i$ and $y_i$ equal to either $\max(r\lambda_i,p_i)$ or $\min(r\lambda_i,p_i)$ for a positive integer parameter $r$. Examples of such tropical coordinate features include:
\begin{alignat*}{4}
    & F_1=\max_i\,\lambda_i \qquad &&  
    F_2 =\max_{i<j}(\lambda_i+\lambda_j) \\
    & F_3=\max_{i<j<k}(\lambda_i+\lambda_j + \lambda_k) \qquad &&  F_4=\max_{i<j<k<l}(\lambda_i+\lambda_j+\lambda_k+\lambda_l) \\
    &F_5=\sum_i \lambda_i \qquad && F_6 =\sum_i \min(r\lambda_i,p_i),
\end{alignat*}
along with the somewhat more complicated
\[ F_7=\sum_j \left[\max_i\big(\min(r \lambda_i,p_i)+\lambda_i\big)-(\min(r \lambda_j,p_j)+\lambda_j)\right].
\]
These seven tropical coordinates were used in \cite{Kalisnik2019} for performing classification on the MNIST database, with $r=28$. 

\noindent\sepline

The third and final algebraic vectorization considered here is generated by extracting complex polynomials from barcodes \cite{Ferri1999, DiFabio2015}. In what follows, the symbol $i$ should be interpreted as $\sqrt{-1}$ (and not as an index for the intervals in $B$). Consider the three continuous maps $R,S,T:\R^2 \to \mathbb{C}$ defined as follows:
\begin{align*}
R(x,y) &= x+iy \\
S(x,y) &= \begin{cases} \frac{y-x}{\alpha\sqrt{2}}\cdot(x+iy)& \text{if }(x,y) \neq (0,0) \\ 0 & \text{otherwise} \end{cases} \\ 
T(x,y) &= \frac{y-x}{2} \cdot \Big[(\cos \alpha - \sin \alpha) + i(\cos \alpha + \sin \alpha)\Big],
\end{align*}
where $\alpha$ is the norm $\sqrt{x^2+y^2}$.

\begin{definition}\label{def:compoly}
Given a barcode $\mu:B \to \Z_{>0}$, let $X:\R^2 \to \mathbb{C}$ be any one of the three functions $R,S,T$ defined above. The {\bf complex polynomial} vectorization of $\mu$ of type $X$ is the sequence of coefficients of the complex polynomial in one variable $z$ given by
\[
C_X(z) := \prod_{[p,q] \in B} \left[z-X(p,q)\right]^{\mu_{p,q}}.
\]
\end{definition}
\noindent In practice, it is customary to either take only the first few highest degree coefficients of $C_X(z)$ or to multiply it by a suitable power of $z$. This is done to guarantee that the feature vectors assigned to a collection of barcodes all have the same dimension.

{\bf Other Algebraic Vectorizations:} In the subsequent section, we describe how to extract vectorizations by using barcode data to build curves which take values in a vector space. Once such a curve has been extracted, one can compute its {\em path signature} via iterated integrals \cite{chevkor}. The path signature resides in the tensor algebra of the target vector space; elements of the tensor algebra are equivalent to coefficients of non-commuting polynomials, and hence constitute algebraic vectorizations of barcodes --- see  \cite{Nanda2020, giusti2021} for examples of this approach.

\subsection{Curve Vectorizations} There are several interesting ways of turning barcodes into one or more  curves, which for our purposes here mean (piecewise) continuous maps from $\R$ to a convenient vector space. Feature vectors can then be constructed by sampling the given curve at finite subsets of $\R$. Perhaps the simplest and most widely used curve-based vectorization is the following.

\begin{definition}\label{def:betticurve}
The {\bf Betti curve} of $\mu:B \to \Z_{>0}$ is the curve $\beta_\mu:\R \to \R$ given by
\[
\beta_\mu(t)=\sum_{[p,q] \in B} \mathbb{1}_{p \leq t < q} \cdot \mu_{p,q}.
\]
\end{definition}

\noindent Here $\mathbb{1}_\bullet$ is the indicator function as described in Definition \ref{def:entropy_summary}, so this function counts the number of intervals (with multiplicity) in $B$ which contain $t$. Very similar in spirit (and formula) to the Betti curve is the following vectorization from \cite{Lawson2022}.

\begin{definition}
The {\bf lifespan curve} of $\mu:B \to \Z_{>0}$ is the map $L_\mu:\R \to \R$ given by
\[
L_\mu(t)=\sum_{[p,q] \in B} \mathbb{1}_{p \leq t < q} \cdot \mu_{p,q} \cdot (q-p).
\]
\end{definition}

It is not  difficult to create very different-looking Betti and lifespan curves from two barcodes which have arbitrarily small bottleneck distance --- we can always add lots of very small intervals to a given barcode without changing its bottleneck distance by a significant amount. One way to rectify the bottleneck instability of Betti and lifespan curves is to test the containment not only of $t$ in each interval $[p,q] \in B$, but rather of the largest subinterval of the form $[t-s,t+s]$. This modification leads to one of the oldest and best-known stable curve vectorizations \cite{Bubenik2015,Bubenik2017}, as defined below.

\begin{definition}\label{def:landscape}
The {\bf persistence landscape} of the barcode $\mu:B \to \Z_{>0}$ is a sequence of curves $\set{\Lambda^\mu_i:\R \to [-\infty,\infty] \mid i \in \Z_{>0}}$ given by
\[
\Lambda^\mu_i(t) := \sup\set{s \geq 0 ~ \Big{|} ~ \left(\sum_{[p,q] \in B} \mathbb{1}_{[t-s,t+s] \subset [p,q]} \cdot \mu_{p,q}\right) \geq i}.
\]
\end{definition}

By convention, the supremum over the empty set is zero. Moreover, since our barcode $B$ is assumed to be finite, the landscape functions $\Lambda^\mu_i$ become identically zero for sufficiently large $i$. An alternate approach to defining persistence landscapes comes from the function $\Delta:B \times \R \to \R$, given by
\begin{align}\label{eq:tent}
\Delta([p,q],t) := \max\left(\min(t-p,q-t),0\right).
\end{align}
For each $i \in \Z_{>0}$, the curve $\Lambda_i^\mu$ from Definition \ref{def:landscape} equals the $i$-th largest number in the multiset that contains $\mu_{p,q}$ copies of $\Delta([p,q],t)$ for each interval $[p,q]$ in $B$. The fourth and final curve vectorization that we consider here was introduced in \cite{Chazal2014}, and it is also defined in terms of the functions $\Delta$ from \eqref{eq:tent}. 
\begin{definition}\label{def:silhou}
Let $w:B \to \R_{>0}$ be any function, which we will denote $[p,q] \mapsto w_{p,q}$. The $w$-weighted {\bf persistence silhouette} of $\mu:B \to \Z_{>0}$ is the map $\phi^w_\mu:\R \to \R$ defined as the weighted average
\[
\phi^w_\mu(t) := \frac{\sum w_{p,q} \cdot \mu_{p,q} \cdot \Delta([p,q],t)}{\sum w_{p,q} \cdot \mu_{p,q}}.
\]
Here both sums on the right are indexed over all $[p,q] \in B$, and $\Delta$ is defined in \eqref{eq:tent}.
\end{definition}

Reasonable choices of weight functions are provided by setting $w_{p,q}=(q-p)^\alpha$ for a real-valued scale parameter $\alpha \geq 0$. For small $\alpha$, the shorter intervals dominate the value of the silhouette curve, whereas for large $\alpha$ it is the longer intervals which play a more substantial role --- see \cite[Sec 4]{Chazal2014} for details.

{\bf Other Curve Vectorizations:} See the {\em envelope embedding} from \cite{Nanda2020}, the {\em accumulated persistence function} in \cite{biscio2019}, and the {\em persistent Betti function} of \cite{xia2018}. In \cite{dong2020}, the persistent Betti function is decomposed along the Haar basis to produce a vectorization. More recently, \cite{Lawson2022} provides a general framework for constructing several different curve vectorizations.

\subsection{Functional Vectorizations}

Here we catalogue those barcode vectorizations which are given by maps from spaces other than $\R$. The first, and perhaps most prominent member of this category is the following vectorization from \cite{Adams2017}. Its definition below makes use of two auxiliary components besides the given barcode $\mu:B \to \Z_{>0}$. The first is a continuous, piecewise-differentiable function $f:\R^2 \to \R_{\geq 0}$ satisfying $f(x,0)=0$ for all $x \in \R$. And the second is a collection of smooth probability distributions $\Psi := \set{\psi_{p,q} \mid [p,q] \in B}$ where $\psi_{p,q}$ has mean $(p,q-p)$.

\begin{definition}\label{def:persimages}
The {persistence surface} of $\mu:B \to \Z_{> 0}$ with respect to $f$ and $\Psi$ (as described above) is the function $\R^2 \to \R$  given by
\[
\rho^\mu_{f,\Psi}(x,y)=\sum_{[p,q] \in B} \mu_{p,q} \cdot f(p,q-p) \cdot \psi_{p,q}(x,y).
\]
The {\bf persistence image} $\mathbf{I}^\mu_{f,\Phi}$ of $\mu$ with respect to $(f,\Phi)$ assigns a real number to every subset $Z \subset \R^2$; this number is given by integrating the persistence surface over $Z$:
\[
\mathbf{I}^\mu_{f,\Psi}(Z)=\iint_Z ~ \rho^\mu_{f,\Psi}(x,y) ~ dx~dy.
\]
\end{definition}

In order to obtain a vector from the persistence image, one lets $Z$ range over grid pixels in a rectangular subset of $\R^2$ and renormalizes the resulting array of numbers, thus producing a grayscale image. Standard choices of $f$ and $\Psi=\set{\psi_{p,q}}$ are:
\begin{align*}
f(x,y) &= \begin{cases} 0 & t \leq 0 \\ t/\lambda_\text{max} & 0 < t < \lambda_\text{max} \\ 1 & t > \lambda_\text{max}\end{cases} \\
\psi_{p,q}(x,y) &= \frac{1}{2\pi\sigma^2} \cdot  \exp\left(-\frac{(x-p)^2+(y-(q-p))^2}{2\sigma^2}\right).
\end{align*}
Here $\lambda_\text{max}$ is the largest lifespan $\max_{[p,q] \in B} (q-p)$ encountered among the intervals in $B$, and $\sigma$ is a user-defined parameter which forms the common standard deviation of every $\psi_{p,q}$ in sight.

The second and final functional vectorization which we will examine was introduced in the paper \cite{Perea2022}. Set $\mathbb{W} := \set{(x,y) \in \R^2 \mid 0 \leq x < y}$, and note that points $(x,y) \in \mathbb{W}$ parameterize intervals $[x,y] \subset \R$ with strictly positive length that could possibly lie in a given barcode. Let $C_c(\mathbb{W})$ be the set of all continuous functions $f:\mathbb{W} \to \R$ with compact support\footnote{In other words, $C_c(\mathbb{W})$ contains those continuous real-valued functions on $\mathbb{W}$ which evaluate to $0$ outside the intersection of a sufficiently large rectangle with $\mathbb{W}$ in $\R^2$.}. The given barcode $\mu:B \to \Z_{>0}$ induces a function $V_\mu:C_c(\mathbb{W}) \to \R$ via
\begin{align}\label{eq:tempmap}
V_\mu(f)=\sum_{[p,q] \in B} \mu_{p,q} \cdot f(p,q-p).
\end{align}
A subset $T$ of $C_c(\mathbb{W})$ is called a {\em template system} if for any distinct pair $\mu_1:B_1 \to \Z_{>0}$ and $\mu_2: B_2 \to \Z_{>0}$ of barcodes, there exists at least one $f \in T$ so that $V_{\mu_1}(f) \neq V_{\mu_2}(f)$. 

\begin{definition}\label{def:tempfunc}
Fix an integer $n > 0$ and let $\text{\rm Sub}_n(T)$ be the collection of all size $n$ subsets of a template system $T$ as described above. The {\bf template function} vectorization of $\mu:B \to \Z_{>0}$ with respect to $T$ is the map $\tau:\text{\rm Sub}_n(T) \to \R^n$ defined as follows. Given $f=\set{f_1,\ldots,f_n}$ in $\text{\rm Sub}_n(T)$, the associated vector in $\R^n$ is 
\[
\tau^\mu(f) := \left(V_\mu(f_1),\ldots,V_\mu(f_n)\right), 
\]
where $V_\mu(f_i)$ is as defined in \eqref{eq:tempmap}.
\end{definition}
Two convenient choices of $T$, called {\em tent functions} and {\em interpolating polynomials}, have been highlighted in \cite{Perea2022}. Tent functions are indexed by points $(u,v) \in \R^2$ and require an additional parameter $\delta > 0$; they have the form
\begin{align}\label{eq:tent2}
g^\delta_{u,v}(x,y)=\max\left(1-\frac{1}{\delta} \cdot \max(|x-u|,|y-v|),0\right)
\end{align}
By construction, each such function is supported on the square of side length $2\delta$ around the point $(u,v)$ in the birth-lifespan plane. The normal pipeline for selecting finitely many template functions requires covering a sufficiently large bounded subset of $\mathbb{W}$ with a square grid and then selecting the appropriate tent functions supported on grid cells. We direct interested readers to \cite[Sections 6 and 7]{Perea2022} for details on interpolating polynomials and for suggestions on how one might select suitable $n$ and $f \in \text{\rm Sub}_n(T)$ for a given classification task.

{\bf Other Functional Vectorizations:} See the {\em generalised persistence landscape} in \cite{berry2020} and the {\em crocker stacks} of \cite{xian2020}. 

\subsection{Ensemble Vectorizations} Our last category contains two methods which require access to a sufficiently large collection of {\em training} barcodes $\mu_i:B_i \to Z_{>0}$ in order to generate a vectorization. The first of these methods, introduced in \cite{polanco2019adaptive}, is a modification of the template system vectorization from Definition \ref{def:tempfunc}. We recall that $\mathbb{W} \subset \R^2$ is defined as $\set{(x,y) \mid 0 \leq x < y}$ and that every barcode $B$ is identified with a subset $P(B) \subset \mathbb{W}$ via the map that sends intervals $[p,q]$ of positive length to points $(p,q)$.

\begin{definition}\label{def:adaptemp} The {\bf adaptive template system} induced by a collection of barcodes $\set{\mu_i:B_i \to \Z_{>0}}$ is obtained via the following two steps. Letting $P \subset \mathbb{W}$ be the union $\bigcup_i P(B_i)$, one
\begin{enumerate}
    \item identifies finitely many ellipses $E_j \subset \mathbb{W}$ which tightly contain $P$, and then
    \item constructs suitable functions $g_j$ supported on $E_j$, as described in \eqref{eq:ellfunc} below.
\end{enumerate}
\end{definition}

The desired vectorization of a new barcode $\mu:B \to \Z_{>0}$ is now obtained by using these $g_j$, rather than tent functions, as template functions in Definition \ref{def:tempfunc}. Three different methods for finding the $E_j$ can be found in \cite[Sec 3]{polanco2019adaptive}. Let $v^*$ denote the transpose of a given vector $v$ in $\R^2$. Now each ellipse $E$ with centre $x = (x_1,x_2)^*$ corresponds to a symmetric $2 \times 2$ matrix $A$ satisfying
\[
E = \set{z \in \R^2 \mid (z-x)^*A(z-x) = 1}.
\]
Setting $h(z) := (z-x)^*A(z-x)$, the adaptive template function $g$ supported on $E$ is
\begin{align}\label{eq:ellfunc}
    g(z) = \begin{cases}
    1-h(z) & h(z) < 1 \\
    0 & \text{otherwise}.
\end{cases}     
\end{align}

The second instance of an ensemble vectorization framework which we benchmark in this paper is from \cite{Royer2021}. Let $\mu_i:B_i \to \Z_{>0}$ be a collection of training barcodes as before, and fix a dimension parameter $b \in \Z_{>0}$. Much like the adaptive template systems of Definition \ref{def:adaptemp}, the {\em automatic topology-oriented learning} (ATOL) vectorization is a two-step process for mapping each $B_i$ to a vector space, which in this instance is always $\R^b$.

\begin{definition}\label{def:atol}
The {\bf ATOL contrast functions} corresponding to the collection of barcodes $\set{\mu_i:B_i \to \Z_{>0}}$ and parameter $b \in \Z_{>0}$ are obtained as follows:
\begin{enumerate}
    \item Treating the point clouds \[P_i := \set{(p,q) \in \R^2 \mid [p,q] \in B_i \text{ and }q>p}\] as discrete measures on $\R^2$, one estimates their average measure $E$.
    \item Let ${\bf z} := (z_1,z_2,\ldots,z_b)$ be a point sample in $\R^2$ drawn (in independent, identically distributed function) along $E$. Define the real numbers $\sigma_i({\bf z})$ for $1 \leq i \leq b$ by
    \[
    \sigma_i({\bf z}) := \frac{1}{2} \max_{j \neq i} \|z_j-z_i\|_2,
    \]
    where $\|\bullet\|_2$ denotes the usual Euclidean norm on $\R^2$.
\end{enumerate}
The contrast functions $\set{\Omega_i:\R^2 \to \R \mid 1 \leq i \leq b}$ are now given by
\[
\Omega_i(x) = \exp\left(-\frac{\|x-z_i\|}{\sigma_i({\bf z})}\right).
\]
\end{definition}
The reader is directed to \cite[Algorithm 1]{Royer2021} for further details. Once the contrast functions have been produced in the manner described above, the corresponding {\bf ATOL vectorization} of a given barcode $\mu:B \to \Z_{>0}$ equals $\left(\Omega_1^\mu,\ldots,\Omega_b^\mu\right)$, where
\[
\Omega_i^\mu := \sum_{[p,q] \in B} \mu_{p,q} \cdot \Omega_i(p,q).
\]

{\bf Other Ensemble Vectorizations: } The {\em persistence codebooks} approach from \cite{Zielinski2021} proposes three different types of barcode vectorizations; these are based on  bag-of-word embeddings, VLAD (vector of locally aggregated descriptors), and Fisher Vectors respectively.

\section{Datasets}\label{sec:datasets}

The vectorization methods described in the preceding section have been benchmarked against three standard datasets; these are described below and arranged in increasing order of difficulty for topological methods. All three of them have been used in the past for comparing vectorizations (or kernels)  for persistence barcodes \cite{Perea2022, polanco2019adaptive, reininghaus2015, carriere2017sliced, guo2018sparse, Nanda2020}.

\begin{figure}[h!]
    \includegraphics[width=.75\textwidth]{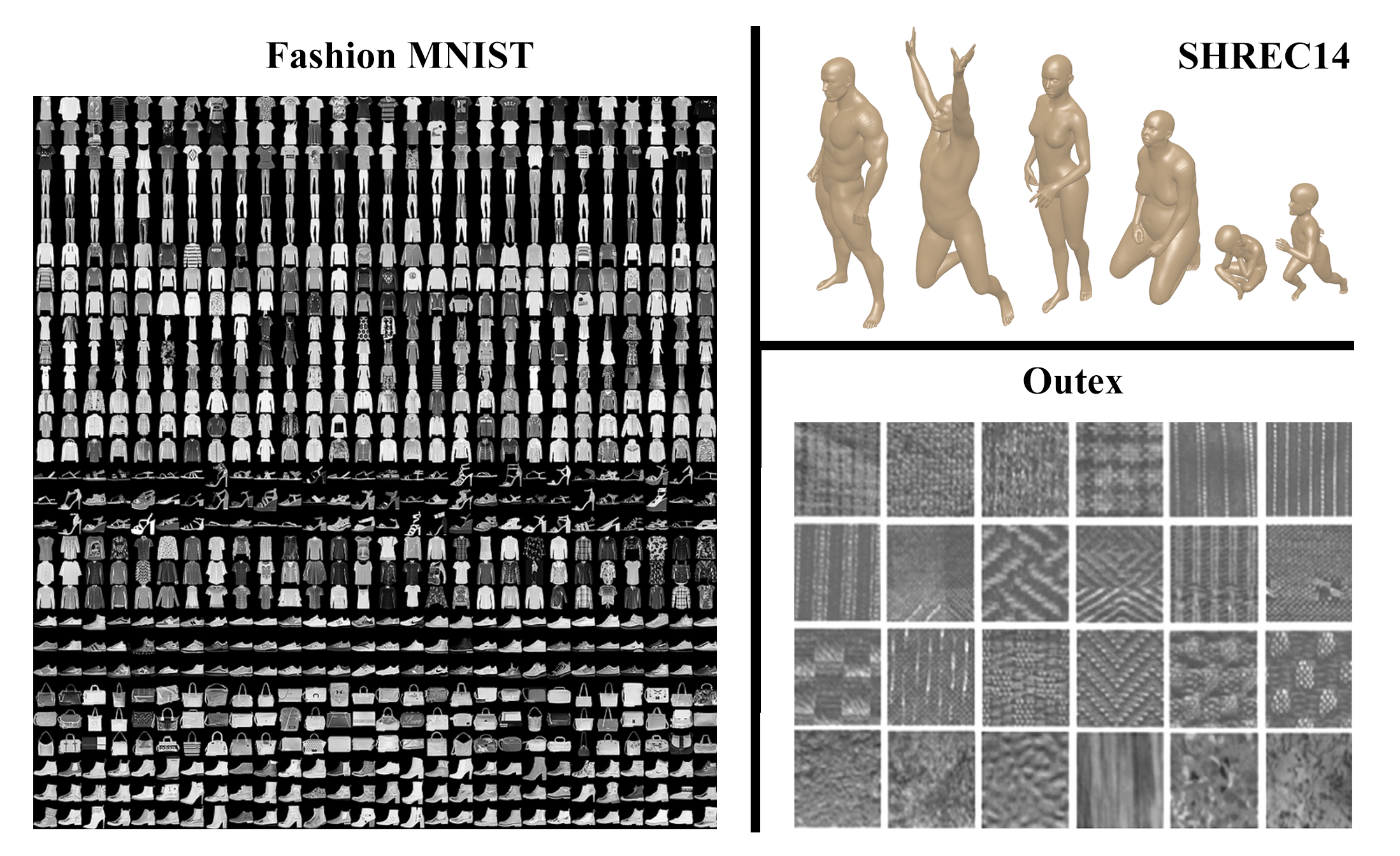}
    \caption{Samples from datasets used in our experiments}\label{datasets}
    \label{fig:galaxy}
\end{figure}

\subsection{Outex}
Outex is a database of images developed for the assessment of texture classification algorithms \cite{ojala2002outex} --- see Fig. \ref{datasets}, right-bottom, for some samples of textures from the 68 categories. Each texture class contains $20$ images of size $128 \times 128$ pixels, which results in $1,360$ images in total. We designed a reduced version of the experiment by randomly selecting $10$ of the total $68$ classes in the dataset, which we refer to as \textbf{Outex10} below. The full classification is referred to as \textbf{Outex68}. In both cases, a train/test split of 70/30 has been applied. 

We treat each image as a cubical complex; the filtration is induced by considering the pixel intensity on the 2-dimensional cells, which is inherited by other cells via the lower-star and upper-star filtrations. Persistent homology barcodes are computed in dimensions $0$ and $1$ using the GUDHI library \cite{gudhi:CubicalComplex}. No pre-processing has been applied to the images.

\subsection{SHREC14}

The Shape Retrieval of a non-rigid 3D Human Models dataset, usually abbreviated SHREC14 \cite{Pickup2014}, is designed to test shape classification and retrieval algorithms. It contains real and synthetic human shapes and poses stored as 3D meshes (which are already simplicial complexes). We use the synthetic part of the dataset; this constitutes a classification task with 15 classes ($5$ men, $5$ women and $5$ children), each one with $20$ different poses --- see the upper-right corner of Fig. \ref{datasets}.

 We apply the Heat Kernel Signature (HKS) to obtain filtrations  \cite{sun2009HKS, reininghaus2015}. For a fixed real parameter $t > 0$, this filtration assigns to each mesh point $x$ the value
 \begin{align}\label{eq:hks}
 \text{HKS}_t(x) = \sum_{i=0}^\infty e^{-\lambda_it} \cdot \phi_i(x)^2
 \end{align}
 Here $\lambda_i$ and $\phi_i$ are eigenvalues and corresponding eigenfunctions of (a discrete approximation to) the Laplace-Beltrami operator of the given mesh. Every simplex of dimension $> 0$ is assigned the largest value of $\text{HKS}_t$ encountered among its vertices. We used the pre-computed barcodes (for such filtrations across a range of $t$-values) which have been provided in the repository\footnote{\url{https://github.com/barnesd8/machine_learning_for_persistence}} accompanying \cite{Perea2021}. Of the $300$ samples, $70\%$ were used for training and the other $30\%$ for testing.

\subsection{FMNIST}

The Fashion-MNIST database contains $28\times 28$ grayscale images ($7,000$ images per class, with $10$ classes) --- see the left side of Fig. \ref{datasets} for some sample images. We split this dataset into $60,000$ training and $10,000$ testing images.

The filtration used for generating barcodes is as follows: we performed  padding, median filter, and shallow thresholding before computing {\em canny edges} \cite{canny}. Then each pixel is given a filtration value equalling its distance from the edge-pixels. Finally, all other cells inherit filtration values from the top pixels via the lower star filtration rule.

\section{Results}  \label{sec:results}                       
Here we report the classification accuracy of the thirteen vectorization methods from Section \ref{sec:vecs} on each of the three datasets from Section \ref{sec:datasets}. Implementation details and parameter choices are provided in Appendix \ref{sec:impl}. The source code is available at the following GitHub repository: \url{https://github.com/Cimagroup/vectorization-maps}.

\subsection{Outex}

Table \ref{tab:outex10} displays the classification accuracy for the smaller (and easier) experiment on $10$ classes. As one might expect, all techniques perform rather well, with {\bf Persistence Statistics} and {\bf Algebraic Functions} sharing the best performance with $99.2\%$ accuracy each, followed closely by Persistent Silhouettes with $98.3\%$ each. 

Results from the full experiment with $68$ classes are contained in Table \ref{tab:outex68}; as one might expect, the performance of every single vectorization degrades in the passage from Outex10 to Outex68. Here {\bf Persistence Statistics} is the clear winner by a significant margin, earning $93.4\%$ accuracy. Tropical Coordinates ranks second with $88.7\%$. Setting aside the outstanding performance of Persistence Statistics, it appears clear from these results that the algebraic vectorizations perform far better on Outex68 than the vectorizations from the other categories.

\begin{table}[ht]
    \centering
    \begin{tabular}{|l|l|l|l|}
    \hline
       {\bf  Vectorization Method }&{\bf  Accuracy} & {\bf Parameters} & {\bf Estimator}\\ \hline\hline\hline
        Persistence Statistics & {\bf 0.992} & ~ & SVM, rbf kernel, $C_1$, $\gamma_1$
        \\ \hline 
        Entropy Summary& 0.975 & 100 &  SVM, rbf kernel,  $C_1$, $\gamma_1$ \\\hline \hline
        Algebraic Functions &{\bf 0.992} & ~ &  SVM, linear kernel, $C_3$ \\ \hline
        Tropical Coordinates & 0.975 & 250 & SVM, linear kernel, $C_4$ \\ \hline
        Complex Polynomial & 0.950 &  5,\,R &  SVM, rbf kernel,  $C_1$, $\gamma_1$  \\ \hline \hline
        Betti Curve & 0.908 & 200 &  SVM, rbf kernel,  $C_1$, $\gamma_1$  \\ \hline
        Lifespan Curve & 0.975 & 100 &  SVM, rbf kernel,  $C_1$, $\gamma_1$ \\ \hline
        Persistence Landscape & 0.975 &  50,\,20 &  SVM, rbf kernel, $ C_2$, $\gamma_2$ \\ \hline
        Persistence Silhouette & 0.983 &  100, 0 &  SVM, rbf kernel,  $C_1$, $\gamma_1$  \\ \hline\hline
        Persistence Image & 0.938 & 1, 25 & RF,   n=500 \\ \hline
        Template Function & 0.958 &  35,\,20 &  SVM, rbf kernel,  $C_1$, $\gamma_1$  \\ \hline\hline
        Adaptive Template System & 0.975 &  GMM, 40 &  SVM, rbf kernel,  $C_1$, $\gamma_1$  \\ \hline
         ATOL & 0.967 & 32 &  SVM, linear kernel, $C_4$ \\ \hline
    \end{tabular}
    \vspace{0.2cm}
    \caption{  {\bf Outex10 results}. The relevant parameter values are $C_1=936.5391$, $\gamma_1=0.0187$, $ C_2=914.9620$, $\gamma_2=0.0061$, $C_3=86.0442$, and $C_4=998.1848$.}\label{tab:outex10}
\end{table}
\begin{table}[!h]
    \centering
    \begin{tabular}{|l|l|l|l|}
    \hline
        {\bf Vectorization Method} & {\bf Accuracy} & {\bf Parameters} & {\bf Estimator} \\ \hline\hline\hline
        Persistence Statistics & {\bf 0.934} & ~ &  SVM, rbf kernel, $C_1$, $\gamma_1$
        \\ \hline
        Entropy Summary & 0.859 & 100 &  SVM, poly kernel, $C_2$, $\gamma_2$, deg=2 \\ \hline\hline
        Algebraic Functions & 0.875 & ~ &  SVM, linear kernel, $C_4$\\ \hline
        Tropical Coordinates & 0.887 & 50 &  SVM, linear kernel, $C_5$\\ \hline
        Complex Polynomial & 0.846 &  10, R & SVM, linear kernel, $C_4$ \\ \hline\hline
        Betti Curve & 0.804 & 200 &  SVM, rbf kernel, $C_1$, $\gamma_1$\\ \hline
        Lifespan Curve & 0.842 & 100 &  SVM, rbf kernel, $C_1$, $\gamma_1$ \\ \hline
        Persistence Landscape & 0.822 &  50, 20 & SVM, rbf kernel, $C_3$, $\gamma_3$  \\ \hline
        Persistence Silhouette & 0.844 &  100, 1 & SVM, linear kernel, $C_4$ \\ \hline\hline
        Persistence Image & 0.762 & 1, 150 & RF, n=500 \\ \hline
        Template Function & 0.831 & 35, 20 & RF, n=200  \\ \hline\hline
        Adaptive Template Sys. & 0.819 & GMM, 50 &  SVM, linear kernel, $C_6$ \\ \hline
        ATOL & 0.854 & 16 &  SVM, linear kernel, $C_7$  \\ \hline
    \end{tabular}
      \vspace{0.2cm}
      \caption{  {\bf Outex68 results}. The optimal parameter values are
       $C_1=936.5391$, $\gamma_1=0.0187$, $C_2=957.5357$, $\gamma_2=0.0120 $, $C_3=914.9620$, $\gamma_3=0.0061$, $C_4=998.1848$, $C_5=884.1255$,  $C_6=143.1201$ and  $C_7=494.0596 $.
      }\label{tab:outex68}
\end{table}

We note that the authors of \cite{Lawson2022} have also used Outex to compare the performance of various curve vectorizations, with Persistence Statistics being used as a baseline. They also obtained their best results with Persistence Statistics.


\subsection{SHREC14}

We used $10$ different $t$-values $t_1 < t_2 < \cdots < t_{10}$, as in \cite{reininghaus2015, Perea2022, polanco2019adaptive}, for generating filtrations via the heat kernel from \eqref{eq:hks}. At $t_{10}$ we found several sparse or empty barcodes, which led us to discard that classification problem. Table \ref{tab:bestshrec14} collects the best performance for each method across the first $9$ values of $t$; it also contains values of the optimal parameters (see Appendix \ref{sec:impl}) and the optimal values of $t$.
\begin{table}[!ht]
    \centering
    \begin{tabular}{|l|l|l|l|}
    \hline
       {\bf Vectorization Method} &{\bf  Accuracy} & {\bf Parameters} & {\bf Estimator} 
        \\ 
        \hline\hline\hline
        Persistence Statistics & {\bf 0.947} &$t_5$ &  RF, n=100\\ \hline
        Entropy Summary & 0.723 &$t_6$, 200 & RF, n=300 \\ \hline\hline
        Algebraic Functions & 0.909 &$t_5$  & RF, n=500 \\ \hline
        Tropical Coordinates & 0.844 &$t_6$, 50 & SVM, linear kernel, $C_5$\\ \hline
        Complex Polynomial & 0.889 &$t_6$, 20, S &  SVM, linear kernel, $C_6$\\ \hline\hline
        Betti Curve & 0.728 &$t_5$, 200 & RF, n=100 \\ \hline
        Lifespan Curve & 0.878 &$t_7$, 200 &  SVM, linear kernel, $C_7$\\ \hline
        Persistence Landscape & 0.889 &$t_6$,  50, 10 &  SVM, rbf kernel, $C_1$, $\gamma_1$ \\ \hline
        Persistence Silhouette & 0.867 &$t_6$, 200, 2 &  SVM, rbf kernel, $C_2$, $\gamma_2$\\ \hline\hline
        Persistence Image & 0.916 &$t_5$, 1, 10 & RF, n=100 \\ \hline
        Template Function &  0.944 &$t_5$,  14, 0.7 &  SVM, rbf kernel, $C_3$, $\gamma_3$ \\ \hline\hline
        Adaptive Template Sys. & 0.889 &$t_5$,  GMM, 15 &  SVM, linear kernel, $C_8$ \\ \hline
        ATOL & 0.933 &$t_8$, 16 & SVM, rbf kernel, $C_4$, $\gamma_4$\\ \hline
    \end{tabular}
        \vspace{0.2cm}
    \caption{Best performance of each method on {\bf SHREC14}. The parameters are $C_1=835.6257$, $\gamma_1=0.0002$, $C_2=212.6281$, $\gamma_2=0.0031$, $C_3=879.1425$, $\gamma_3=0.0010$, $C_4=936.5391$, $\gamma_4=0.0187$, $C_5=141.3869$,  $C_6=625.0300$, $C_7=998.1848$, $C_8=274.500$.
      }\label{tab:bestshrec14}
\end{table}

{\bf Persistence Statistics} yielded the best classification accuracy of $94.7\%$, followed closely by Template Functions at $94.4\%$. One remarkable feature of these results is that the dataset does not appear to favour any one category of vectorizations over the other --- it is possible to achieve over $88\%$ accuracy by using a suitable statistical, algebraic, curve, functional or ensemble vectorization. In fact, only the curve-based vectorizations failed to achieve over $90\%$ accuracy on this dataset. The variation of classification accuracy with the heat kernel parameter $t$ is discussed in Appendix \ref{app:shrec}.

\subsection{FMNIST}

The results of our experiments on FMNIST are recorded in Table \ref{tab:fmnist}. We note that these experiments only used information contained in the $0$-dimensional barcodes and that the SVM classifier was not used. The classification accuracy of all the methods is much lower than the corresponding figures for the two preceding datasets. Once more, the {\bf Persistence Statistics} vectorization takes the top spot with $74.9\%$ and Template Functions are slightly behind at $74.7\%$

\begin{table}[!ht]
    \centering
    \begin{tabular}{|l|l|l|}
    \hline
        {\bf Vectorization Method} &{\bf  Accuracy} &{\bf Parameters}  \\ \hline\hline\hline
        Persistence Statistics & {\bf 0.749} & ~ \\ \hline
        Entropy Summary & 0.696 & 30  \\ \hline\hline
        Algebraic Functions & 0.710 & ~  \\ \hline
        Tropical Coordinates & 0.696 & 10  \\ \hline
        Complex Polynomial & 0.661 & 10, R \\ \hline\hline
        Betti Curve & 0.618 & 50  \\ \hline
        Lifespan Curve & 0.692 & 30  \\ \hline
        Persistence Landscape & 0.694 & 30, 5  \\ \hline
        Persistence Silhouette & 0.670 & 30, 0   \\ \hline\hline
        Persistence Image & 0.698 & 1, 12 \\ \hline
        Template Functions & 0.747 & 10, 2 \\ \hline\hline
        Adaptive Template System & 0.602 & GMM, 5 \\ \hline
        ATOL & 0.730 & 16  \\ \hline
    \end{tabular}
        \vspace{0.2cm}
        \caption{{\bf  FMNIST results.} All the scores have been achieved for Random Forest classifier with $100$ trees.}\label{tab:fmnist}
\end{table}

One rather surprising aspect of these results is the fact that Adaptive Template Systems performed far worse than ordinary Template Functions despite having recourse to $60,000$ training barcodes. We do not have a clear explanation for this phenomenon, particularly in light of a fairly competitive performance from ATOL (which was also exposed to the same training data). 

\section{Web Application} \label{sec:webapp}

In order to illustrate and visualize the vectorization methods described here, we have built an interactive web application that runs on any modern browser; it is available at
\[
\text{\url{https://persistent-homology.streamlit.app/}}
\]
The app has been built in Python using the Streamlit library together and makes use of several existing Python libraries. The sidebar contains options for selecting different types of input data and displays several options for data visualization. One sample image/point-cloud from each of the three datasets used in this paper has been pre-loaded, but the user is free to upload their own data. Specifications, formatting guidelines, and downloading instructions are available in our GitHub repository:
\[
\text{\url{https://github.com/dashtiali/vectorisation-app}}
\]

\begin{figure}[h!]
    \includegraphics[scale=.35]{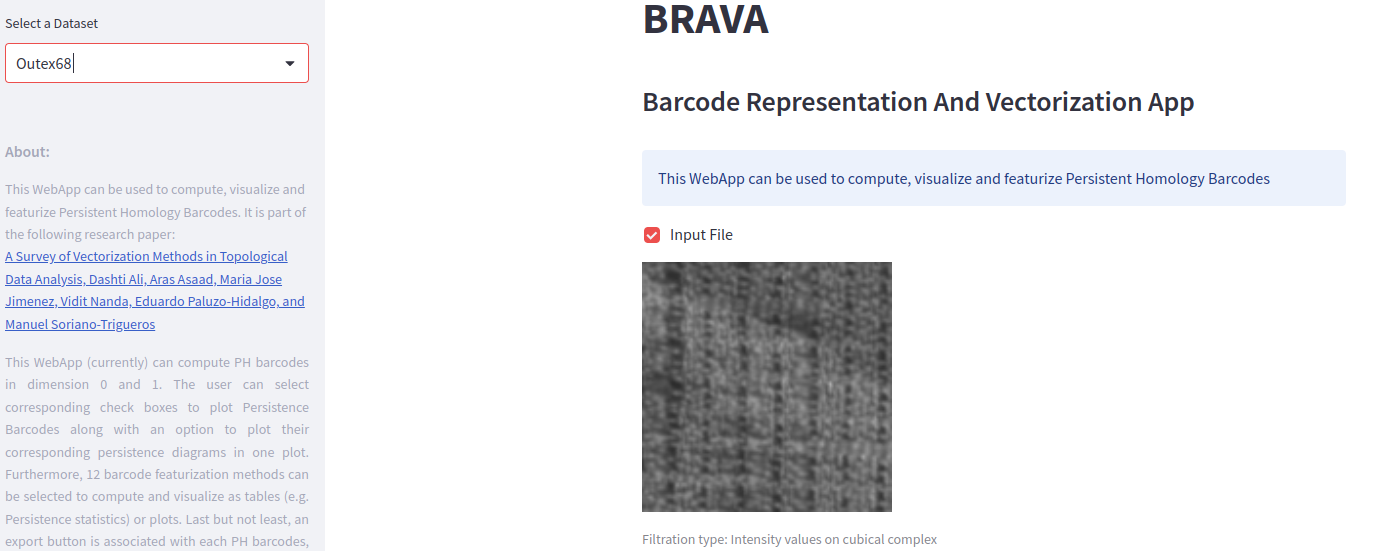}
    \caption{A screenshot of the web app}\label{fig:input}
\end{figure}

\begin{figure}[h!]
    \centering
    \includegraphics[scale=.70]{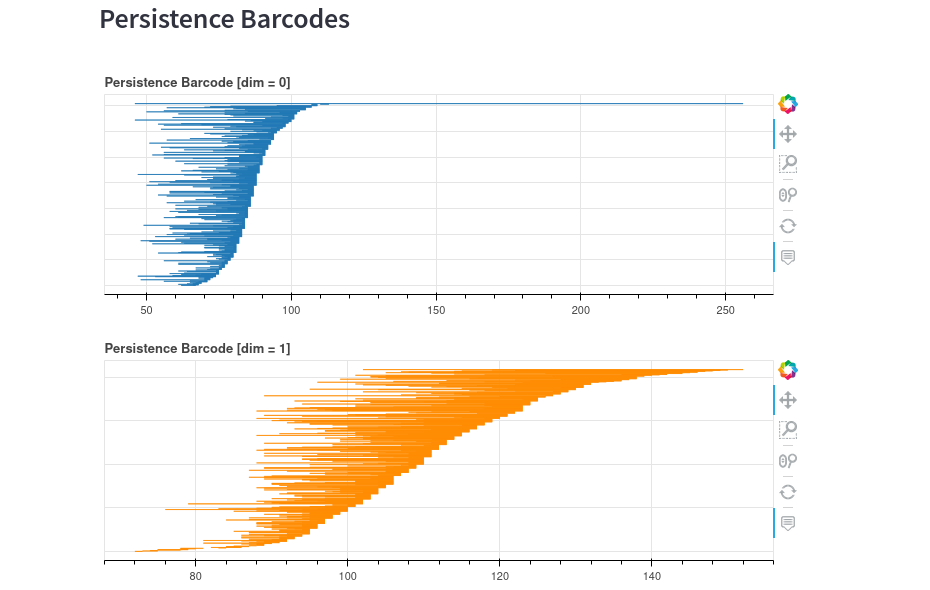}
    \caption{Intervals in barcodes of dimensions $0$ and $1$ as displayed by the web app.}
    \label{fig:barcodes}
\end{figure}

All of the barcode vectorization methods considered in this paper can be computed and visualized in different formats (tables, bar graphs, scatter plots), depending on the type of vectorization being invoked. Barcodes are computed by default in dimensions $0$ and $1$, and  depicted as in Figure \ref{fig:barcodes}.

\begin{figure}[h!]
    \centering
    \includegraphics[scale=.65]{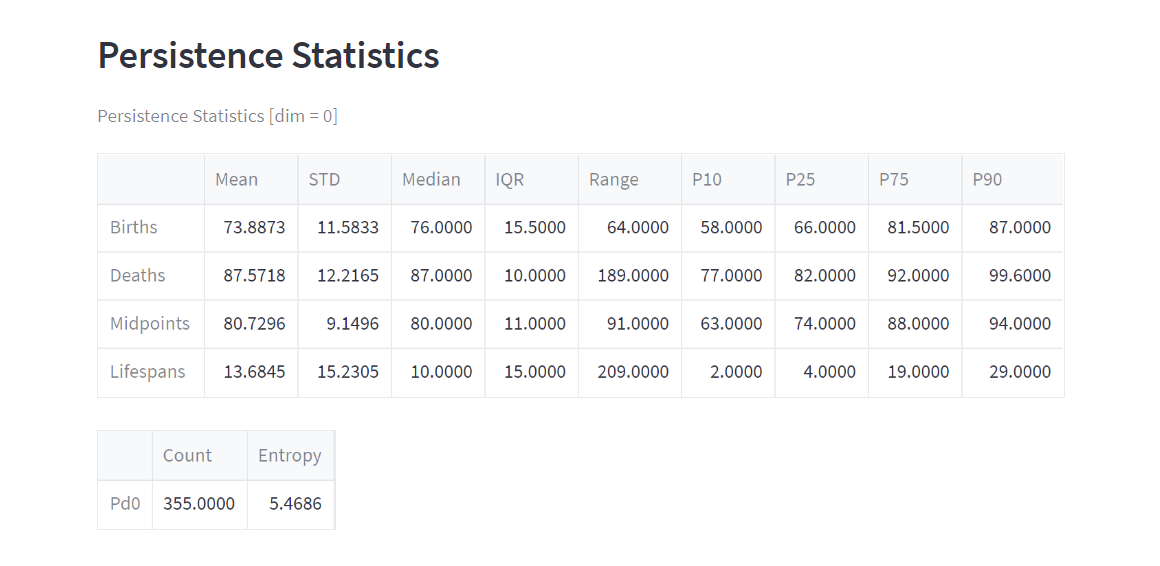}
    \caption{The Persistence Statistics vectorization as shown in the web app.}
    \label{fig:PSt}
\end{figure}

The Persistence Statistics vectorization is purely numerical, so we show its values in a table, as in Figure \ref{fig:PSt}.

\begin{figure}[h!]
    \centering
    \includegraphics[scale=.45]{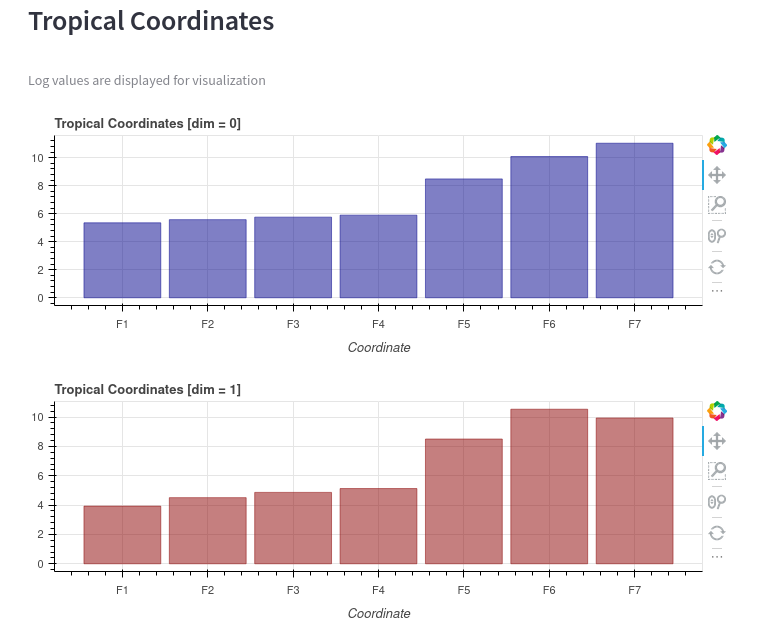}
    \caption{A visualization of the Tropical Coordinates vectorization from the web app.}
    \label{fig:TC}
\end{figure}

Algebraic vectorizations are illustrated as bar graphs. In Figure \ref{fig:TC}, for instance, one finds bars whose heights correspond to values attained by the $7$ chosen tropical coordinate polynomials on the input barcodes.

\begin{figure}[h!]
    \centering
    \includegraphics[scale=.85]{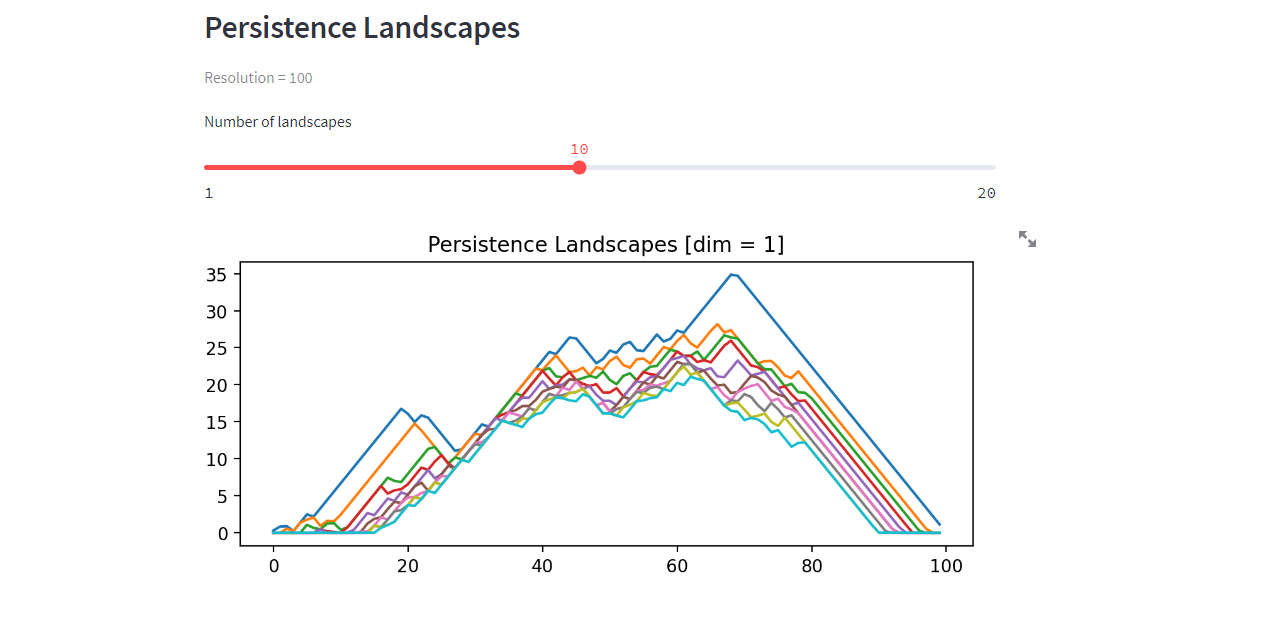}
    \caption{Persistence landscapes in the web app}
    \label{fig:PL}
\end{figure}

Curve vectorizations, such as persistence landscapes, are depicted via piecewise-linear graphs (see Figure \ref{fig:PL}). Sliders have been provided to set the resolution parameter.

\begin{figure}[h!]
    \includegraphics[scale=.55]{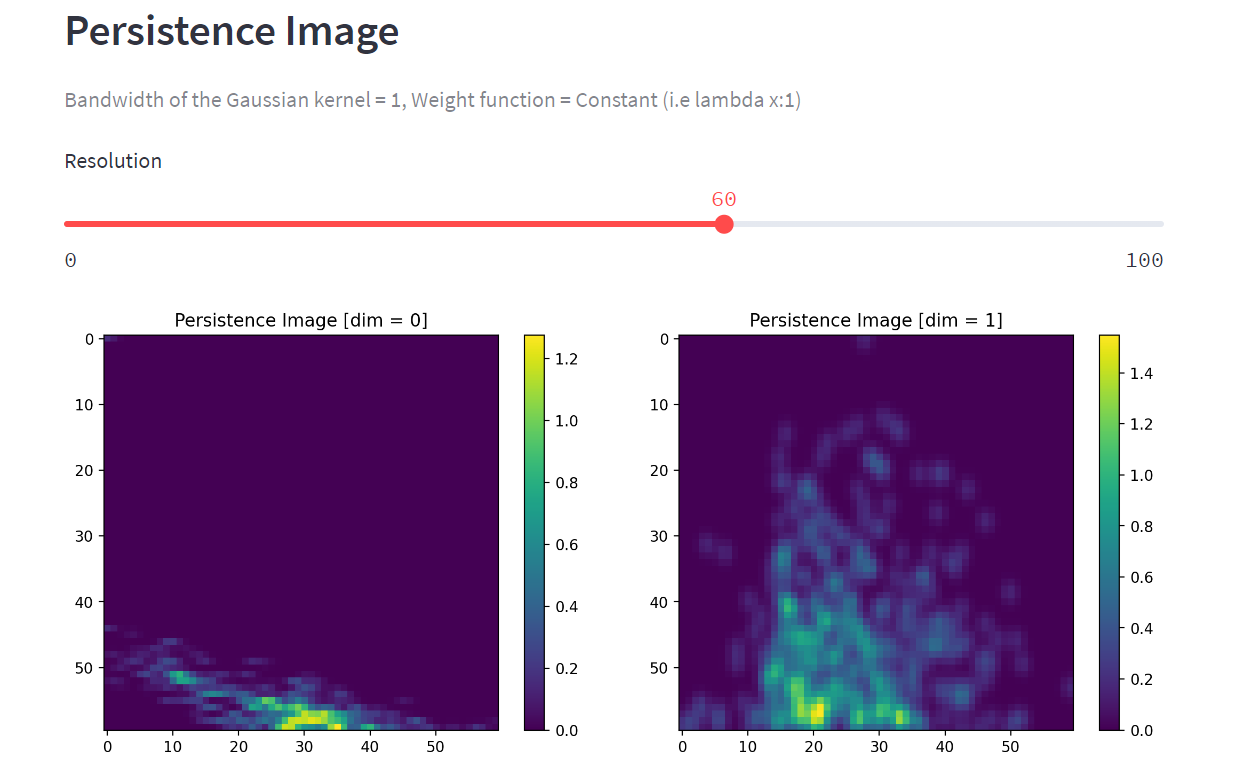}
    \caption{Persistence images as shown in the web app}
    \label{fig:PI}
\end{figure}

Persistence images are displayed as heat maps --- see Figure \ref{fig:PI}.

\begin{figure}[h!]
    \centering
    \includegraphics[scale=.45]{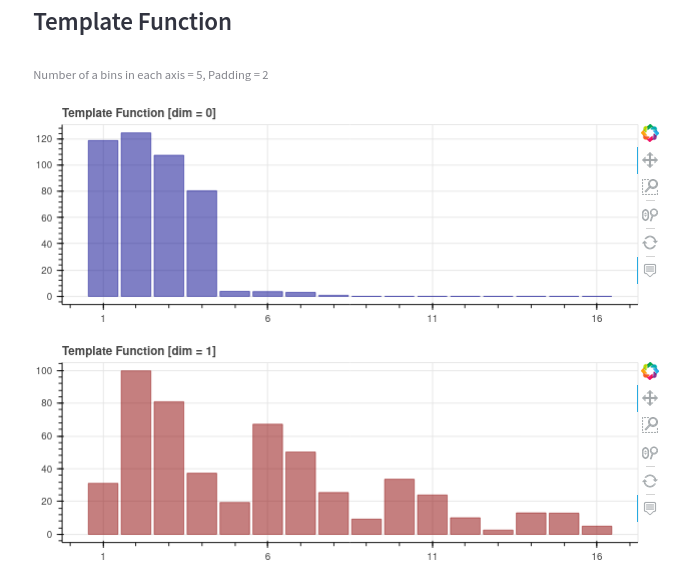}
    \caption{The web app visualization of template functions}
    \label{fig:TF}
\end{figure}

Template Functions, their adaptive version, and ATOL are all displayed as bar graphs with heights of bars indicating the values of the selected functions. Figure \ref{fig:TF}, for instance, depicts Template Functions. 

It is our hope that users will benefit from the ability to generate these visualizations without having to write any code of their own. In order to facilitate downstream analysis, the web app also provides the ability to download the vectors generated by each vectorization method.

\section{Concluding Remarks} \label{sec:conclusion}

At the time of writing, it remains difficult to accurately pinpoint those attributes which might make a given vectorization method a good choice for a particular classification problem. There are no powerful theorems or immutable doctrines available to guide scientists who wish to incorporate topological information into machine learning pipelines. In the absence of such theoretical foundations, the best that one can expect are principled heuristics supported by reproducible empirical evidence. This paper is an outcome of our efforts to provide such evidence. En route, we have organized thirteen available vectorization methods into five categories in Section \ref{sec:vecs} and provided a web application which will allow others to conduct their own experiments involving these methods.

One possible conclusion that may be drawn from the results of Section \ref{sec:results} is that we can dispense with sophisticated vectorization techniques and only use (some variant of) Persistence Statistics. We do not necessarily suggest such a course of action. While it is certainly true that Persistence Statistics earned top honors in all of our experiments and is much faster to compute than the alternatives, there are other factors to consider. In particular, no comparative study such as ours can be truly exhaustive. There is always the chance that making different choices -- for instance, using another dataset for classification, or adding some new polynomials to one of the algebraic vectorizations -- could dramatically update our priors about which methods perform best. 

{\footnotesize
\section*{Acknowledgments}
M.J. Jimenez, E. Paluzo-Hidalgo and M. Soriano-Trigueros are funded by the Spanish grants Ministerio de Ciencia e Innovacion - Agencia Estatal de Investigacion/10.13039/501100011033, PID2019-107339GB-I00 and Agencia Andaluza del Conocimiento, PAIDI-2020 P20-01145. M.J. Jimenez is also funded by a grant of Convocatoria de la Universidad de Sevilla para la recualificacion del sistema universitario español, 2021-23, funded by the European Union, NextGenerationEU.

V. Nanda is supported by EPSRC grant EP/R018472/1 and by
US AFOSR grant FA9550-22-1-0462.

We are grateful to the team of GUDHI and TEASPOON developers, for their work and their support. We are also grateful to Streamlit for providing extra resources to deploy the web app online on Streamlit community cloud.
}

\bibliographystyle{abbrv}
\bibliography{bibliography}

\appendix

\newpage

\section{Implementation and Parameter Details}\label{sec:impl}                                                               
We have made use of several existing software packages, such as GUDHI \cite{gudhi:urm}, Teaspoon\footnote{\url{https://lizliz.github.io/teaspoon/index.html}} or Scikit-learn \cite{scikit-learn}, as well as our own implementations in some cases. Salient information regarding each method has been provided in the list below. Full details can be found in the GitHub  repository accompanying this paper\footnote{\url{https://github.com/Cimagroup/vectorization-maps}}.  

\subsection{Persistence Statistics}

The persistence statistics vectorization from Definition \ref{def:pers_stats} requires no additional parameters. We have implemented this method ourselves.

\subsection{Entropy Summary Function}
We have used the GUDHI implementation of the entropy summary function from Definition \ref{def:entropy_summary}. There is a single resolution parameter which selects the grid points on which the entropy summary function is sampled.

\subsection{Algebraic Functions} The algebraic functions of Definition \ref{def:algfun} are implemented in the Teaspoon package. For reasons which remain unclear to us, this implementation includes a fifth tropical polynomial $f_5= {\max_i \{(q_i - p_i)\}}$ beyond the four ordinary polynomials $f_1,\ldots,f_4$ which were described after Definition \ref{def:algfun}. We do not expect that removing this function will improve the results described below.  

\subsection{Tropical Coordinates}
We have implemented the tropical polynomials $F_1, \ldots, F_7$ described after Definition \ref{def:tropical}. The parameter $r$ has been optimized over the set $\{10,50,250,$ $500,800\}$ for Outex and SHREC14, and over $\set{10,50,250}$ for FMNIST.

\subsection{Complex Polynomials}
We have used the GUDHI implementation of complex polynomials, which have been described in Definition \ref{def:compoly}. We generated the polynomials with respect to all three of the transformations $R, S, T:\R^2 \to \mathbb{C}$. The number of coefficients used was chosen from $\set{5,10,20}$ for Outex and SHREC14 and $\set{3,5,10}$ for FMNIST.

\subsection{Betti Curve}
The Betti curve vectorization from Definition \ref{def:betticurve} has been implemented in GUDHI, and it only requires a resolution parameter. This parameter was chosen from $\set {50,100,200}$ for Outex and SHREC14 and $\set {15, 30, 50}$ for FMNIST.

\subsection{Lifespan Curve}
We implemented the lifespan curve ourselves, with a resolution parameter optimised across the set $\set {50,100,200}$ for Outex and SHREC14 and across the set $\set {15, 30, 50}$ for FMNIST. 

\subsection{Persistence Landscapes}
We have used the GUDHI implementation of persistence landscapes (see Definition \ref{def:landscape}). The are two parameters to consider: the resolution (to identify the grid points where each landscape is sampled) and the total number of landscapes used. The resolution was optimized over $\set {50,100,200}$ for Outex and SHREC14, and over $\set {15, 30, 50}$ for FMNIST; the number of landscapes ranged over $\set {2,5,10,20}$ for Outex and SHREC14 and over $\set {1,2,3,5}$ for FMNIST.

\subsection{Persistence Silhouette}

We have used the GUDHI implementation of persistence silhouettes (see Definition \ref{def:silhou}). The resolution parameter was optimized over $\set{50,100,200}$ for Outex and SHREC14 and over $\set {15, 30, 50}$ for FMNIST; the weight $w$ ranged over $\{0,1,2,5,10,20\}$ for Outex and SHREC14 and $\set{0,1,2,5}$ for FMNIST.

\subsection{Persistence Images}

Persistence images (from Definition \ref{def:persimages}) have been implemented in GUDHI. The resolution parameter $r$, which results in images of size $r \times r$, ranged over $\set{25,75,150}$ for Outex, over $\set {10,20,40}$ for SHREC14, and over $\set {3,6,12,20}$ for FMNIST. Bandwidth values of the Gaussian kernel ($\sigma$ in Definition \ref{def:persimages}) were chosen from $\set{0.05,1}$ for Outex and from $\set{0.05,0.5,1}$ for both, SHREC14 and FMNIST.

\subsection{Template Functions}
We have used code from the repository\footnote{\url{https://github.com/lucho8908/adaptive_template_systems}} provided with the paper \cite{Perea2022} for computing template functions (see Definition \ref{def:tempfunc}). We use tent functions as described in \eqref{eq:tent}, which require two parameters: a grid resolution $\delta$ and a padding parameter $\pi$ (for enlarging the area covered by the square grid). We optimized over 
\begin{itemize}
\item $\delta$ in $\set{35,50,65}$ and $\pi$ in $\set{20,25,30}$ for Outex;
\item $\delta$ in $\set{3,4,\ldots,14,15}$ and $\pi$ in $\{0.5,0.6,\ldots,1.1,1.2\}$ for SHREC14;
\item $\delta$ in $\set{2,3,5,10}$ and $\pi$ in $\set{0.5, 1, 2}$ for FMNIST.
\end{itemize}	

\subsection{Adaptive Template Systems}
The implementation of adaptive template systems (Definition \ref{def:adaptemp}) has also been sourced from the same repository as template functions. We have used the Gaussian mixture model for generating ellipsoidal domains, and require only one parameter: the number of clusters. This has been optimized over
\begin{itemize}
    \item $\set{10,20,30,40,50}$ for Outex,
    \item $\set{5,10,15,20,25,30,35,40,45}$ for SHREC14, and 
    \item $\set{3,4, 5,10,15}$ for FMNIST.
\end{itemize}

\subsection{ATOL}

The ATOL vectorization from Definition \ref{def:atol} has been implemented in GUDHI, and it also requires the number of functions $b$ as a parameter. We have optimized this over $\set {2,4,8,16,32,64}$ for Outex and over $\set {2,4,8,16}$ for both SHREC14 and FMNIST.

\subsection{Dimensions, Classifiers and Hyperparameters}

In the case of Outex, we have concatenated vectors arising from  barcodes of dimensions $0$ and $1$; for SHREC14, the vectors computed from only dimension $1$ barcodes performed better, so the results are only reported for them. Finally, only dimension $0$ barcodes were taken to build vectors for FMNIST. We considered both Support Vector Machine (SVM) \cite{cortes1995support} and Random Forest (RF) \cite{ho1995random} classifiers. Due to convergence issues, only RF has been performed for FMNIST. 

For each parameter of each vectorization method, we accomplished a hyperparameter optimization process based on random search (when optimizing SVM and RF jointly) or grid search (for optimizing RF), with 5-fold cross-validation on the training data, to find the best (hyper)parameters for both the machine learning models and the vectorization methods; then, we assigned to each method the parameters with the best average score among all the 5-fold cross-validation scheme; finally the vectorization methods were evaluated on the test dataset $100$ times, to report the average accuracy.

\section{Heat Kernel Parameter Dependence}\label{app:shrec}  

As mentioned in Section \ref{sec:datasets}, the filtration for {\bf SHREC14} is generated using the Heat Kernel Signature \eqref{eq:hks} which depends on a single parameter $t$. In Table \ref{tab:shrec14} we depict the best classification accuracy of each vectorization method across all $9$ values of $t$ which were used in our experiments. 

\begin{table}[!h]
    \centering
    \begin{tabular}{|c|c|c|c|c|c|c|c|c|c|}
\hline
        Method &$t_1$ &$t_2$ &$t_3$ &$t_4$ &$t_5$ &$t_6$ &$t_7$ &$t_8$ &$t_9$ 
        \\ \hline\hline\hline
        Pers Stat & 0.729 & 0.785 & 0.662 & 0.704 & {\bf 0.947} & 0.910 & {\bf 0.915} &  0.915 & 0.908 \\ \hline
        Ent Sum & 0.378 & 0.333 & 0.522 & 0.536 & 0.656 & 0.723 & 0.633 & 0.656 & 0.530  \\ \hline \hline
        Alg Fun & 0.467 & 0.456 & 0.556 & 0.567 & 0.909 & 0.878 & 0.863 & 0.833 & 0.711 \\ \hline
        Trop Coord & 0.505 & 0.556 & 0.522 & 0.612 & 0.822 & 0.844 & 0.833 & 0.767 & 0.800 \\ \hline
        Com Poly & 0.322 & 0.456 & 0.400 & 0.467 & 0.856 & 0.889 & 0.844 & 0.850 & 0.790 \\ \hline \hline
        Bet Cur & 0.511 & 0.467 & 0.628 & 0.660 & 0.728 & 0.633 & 0.611 & 0.644 & 0.536 \\ \hline
        Lif Cur & 0.456 & 0.411 & 0.593 & 0.639 & 0.789 & 0.833 & 0.878 & 0.833 & 0.798 \\ \hline
        Pers Land & 0.700 & 0.511 & 0.789 & 0.778 & 0.878 & 0.889 & 0.857 & 0.833 & 0.789 \\ \hline
        Pers Sil & 0.400 & 0.378 & 0.556 & 0.589 & 0.811 & 0.867 & 0.856 & 0.856 & 0.656  \\ \hline \hline
        Pers Img & 0.644 & 0.691 & 0.795 & {\bf 0.856} & 0.916 & 0.794 & 0.871 & 0.811 & 0.718 \\ \hline
        Temp Func & 0.778 & 0.735 & {\bf 0.933} & 0.789 & 0.944 & {\bf 0.919} &  0.908 & 0.932 & {\bf 0.922}\\ \hline \hline
        Ad Temp Sys & 0.802 & {\bf 0.872} & 0.833 & 0.727 & 0.889 & 0.856 & 0.889 & 0.844 & 0.633  \\ \hline
        ATOL & {\bf 0.828} & 0.786 & 0.911 & 0.833 & 0.906 & 0.867 & 0.900 & {\bf 0.933} & 0.867 \\ \hline
    \end{tabular}
    \vspace{0.2cm}
    \caption{Best Results for {\bf SHREC14} for various vectorization methods across nine $t$-values.
    }\label{tab:shrec14}
\end{table}

We note that for small values of $t$, the ensemble vectorizations perform best, whereas for intermediate and larger values both ensemble and functional vectorizations achieve good performance. The algebraic and curve based vectorizations perform quite poorly for low $t$-values, but tend to become more competitive between $t_5$ and $t_8$.

\end{document}